\documentclass{ws-ijm}
\usepackage{bbm}
\usepackage{mathrsfs}
\usepackage{xypic}

\newtheorem{thm}{Theorem}[section]
\newtheorem{lem}[thm]{Lemma}
\newtheorem{cor}[thm]{Corollary}
\newtheorem{pro}[thm]{Proposition}
\newtheorem{ex}[thm]{Example}

\newtheorem{defi}[thm]{Definition}

\newcommand{\gl }{{\mathfrak{gl} } }

\newcommand{\lon }{\,\rightarrow\,}
\newcommand{\be }{\begin{eqnarray*}}
\newcommand{\ee }{\end{eqnarray*}}

\newcommand{\defbe}{\triangleq}

\newcommand{\inverse}{^{-1}}

\newcommand{\huaA}{\mathcal{A}}

\newcommand{\huaU}{\mathcal{U}}

\newcommand{\huaW}{\mathcal{W}}

\newcommand{\huaD}{\mathcal{D}}

\newcommand{\huaT}{\mathcal{T}}

\newcommand{\CWM}{C^{\infty}(M)}

\newcommand{\set}[1]{\left\{#1\right\}}

\newcommand{\pibracket}[1]{\left [ #1\right ]_{\pi}}

\newcommand{\frkd}{\mathfrak d}

\newcommand{\frkr}{\mathfrak r}
\newcommand{\frks}{\mathfrak s}

\newcommand{\frky}{\mathfrak y}

\newcommand{\frkL}{\mathfrak L}

\newcommand{\frkX}{\mathfrak X}

\newcommand{\ppairingE}[1]{\left ( #1\right )_E}

\def\gpd{\,\lower1pt\hbox{$\longrightarrow$}\hskip-.24in\raise2pt
         \hbox{$\longrightarrow$}\,}

\newcommand{\connection}{\gamma}
\newcommand{\eomnib}{\mathbf{a}}
\newcommand{\pomnib}{\mathbf{b}}

\newcommand{\LieDerivation}{\frkL}

\newcommand{\pisharp}{\pi}

\newcommand{\half}{\frac{1}{2}}

\newcommand{\Graph}{ \mathbf{G}}

\newcommand{\Rep}{\mathscr{L}}

\newcommand{\conpairing}[1]{\left\langle  #1\right\rangle }

\newcommand{\Dorfman}[1]{\left \{ #1\right \} }

\newcommand{\jet}{\mathfrak{J}}
\newcommand{\jetd}{\mathbbm{d}}
\newcommand{\dev}{\mathfrak{D}}

\newcommand{\EStar}{{E^*}}

\newcommand{\pie}{^\prime}
\newcommand{\Id}{\mathbf{1}}

\newcommand{\e}{\mathbbm{e}}
\newcommand{\p}{\mathbbm{p}}
\newcommand{\id}{\mathbbm{i}}
\newcommand{\jd}{\alpha}

\newcommand{\dM}{\mathrm{d}}

\newcommand{\omni}{\mathcal{E}}
\newcommand{\omnirho}{\rho}

\newcommand{\Hom}{\mathrm{Hom}}
\newcommand{\rank}{\mathrm{rank}}
\newcommand{\Der}{\mathrm{Der}}
\newcommand{\Inn}{\mathrm{Inn}}

\newcommand{\Ext}{\mathrm{Ext}}

\newcommand{\pf}{\noindent{\bf Proof.}\ }
\newcommand{\ad}{\mathrm{ad}}
\newcommand{\Ker}{\mathrm{Ker}}
\newcommand{\Img}{\mathrm{Im}}

\newcommand{\LDorfman}[1]{   \{      #1  \}       }
\newcommand{\falling}[1]{{#1}_{\bullet}}
\begin{document}

\markboth{Z. Chen, Z. Liu and Y. Sheng} {Dirac structures of
omni-Lie algebroids}

\catchline{}{}{}{}{}

\title{DIRAC STRUCTURES OF OMNI-LIE ALGEBROIDS}

\author{ZHUO CHEN}

\address{Department of Mathematics, Tsinghua University, Beijing
100084, China \\
zchen@math.tsinghua.edu.cn}

\author{ZHANGJU LIU}

\address{Department of Mathematics and LMAM, Peking University,
Beijing 100871, China\\
liuzj@pku.edu.cn}

\author{YUNHE SHENG}

\address{Mathematics School $\&$ Institute of Jilin University,
 Changchun 130012, Jilin, China\\
shengyh@jlu.edu.cn}

\maketitle

\begin{abstract}
Omni-Lie algebroids are generalizations of Alan Weinstein's omni-Lie
algebras. A Dirac structure in an omni-Lie algebroid $\dev E\oplus
\jet E$ is necessarily a Lie algebroid together with a
representation on $E$. We study the geometry underlying these Dirac
structures in the light of reduction theory. In particular, we prove
that there is a one-to-one correspondence between reducible Dirac
structures and projective Lie algebroids in $\huaT=TM\oplus E$; we
establish the relation between the normalizer $N_{L}$ of a reducible
Dirac structure $L$ and the derivation algebra $\Der(\pomnib (L))$
of the projective Lie algebroid $\pomnib (L)$; we study the
cohomology group $\mathrm{H}^\bullet(L,\rho_{L})$ and the relation
between $N_{L}$ and $\mathrm{H}^1(L,\rho_{L})$; we describe Lie
bialgebroids using the adjoint representation; we study the
deformation of a Dirac structure $L$, which is related with
$\mathrm{H}^2(L,\rho_{L})$.
\end{abstract}



 \section{Introduction}

Lie algebroids (and local Lie algebras in the sense of Kirillov
\cite{KirillovLocal}) are   generalizations of Lie algebras that
naturally appear in Poisson geometry (and its variations, e.g.,
Jacobi manifolds in the sense of Lichnerowicz
\cite{LichnerowiczJacobi})(see \cite{Mkz:GTGA} for a detailed
description of this subject). Courant algebroids are combinations of
Lie algebroids and quadratic Lie algebras. It was originally
introduced in \cite{CourantDirac} by T. Courant where he first
called them Dirac manifolds, and then were re-named after him in
\cite{LWXmani} (see also an alternate definition
\cite{Uchinoremarks}) by Liu, Weinstein and Xu to describe the
double of a Lie bialgebroid. Recently, several applications of
Courant algebroids and  Dirac structures have been found in
different fields, e.g., Manin pairs and moment maps \cite{AK,Courant
morphism};  generalized complex structures
\cite{BCG,GualtieriGeneralizedComplex}; $L_{\infty}$-algebras and
symplectic supermanifolds \cite{BiSheng,rogersCourant,ROy}; gerbes
\cite{PW} as well as BV algebras and  topological field theories
\cite{BV quantization,ROy1}.

Motivated by an integrability problem of the Courant bracket, A.
Weinstein gives a linearization of the Courant bracket at a point
\cite{Weinomni}, which had been studied from several aspects
recently
\cite{BCG,InteCourant,SLZomilie2algebra,StienonReduction,UchinoOmni}.
 A. Weinstein has shown that an omni-Lie algebra structure
  can encode all Lie algebra structures on a vector space, the next step is,
logically, to find out candidates that could encode all Lie
algebroid structures on a vector bundle.  In a recent work
\cite{CLomni}, we have given a definitive answer to this question.
Over there, a generalized Courant algebroid structure is defined on
the direct sum bundle $\dev E\oplus \jet E$, where $\dev E$ and
$\jet E$ are the gauge Lie algebroid  and the jet bundle of a vector
bundle $E$ respectively. Such a structure is called an {\em omni-Lie
algebroid} since it reduces to the omni-Lie algebra if the base
manifold is a point \cite{Weinomni}. Furthermore, an omni-Lie
algebroid is the first example of $E$-Courant algebroids
\cite{CLSecourant}.

It is well known that the theory of Dirac structures has wide and
deep applications in both mathematics and physics (e.g.,
\cite{BC,port
hamiltonian,Dorfman1993,GualtieriGeneralizedComplex,Hitchin,Wade}).
In \cite{CLomni}, only some special Dirac structures were studied.
The authors proved that there is a one-to-one correspondence between
Dirac structures coming from
   bundle maps
$\jet E \rightarrow \dev E$  and  Lie algebroid (local Lie algebra)
structures on $E$ when $\mathrm{rank}(E)\geq 2$  ($E$ is a line
bundle). In other words, Dirac structures that are graphs of maps
  actually underlines the geometric objects of Lie algebroids, or local Lie algebras.

As a continuation of \cite{CLomni}, the present paper explores what
a general Dirac structure of the omni-Lie algebroid can encode. For
a vector space $V$, Weinstein proved that Dirac structures in the
omni-Lie algebra $\gl(V)\oplus V$ correspond to Lie algebra
structures on subspaces of $V$ in \cite{Weinomni}. For a vector
bundle $E$ over $M$, Dirac structures in the omni-Lie algebroid
$\omni=\dev E\oplus \jet E$ turn out to be more complicated than
that of omni-Lie algebras.  The key concept that we need is a
projective Lie algebroid, which is a subbundle
$A\subset\huaT=TM\oplus E$, equipped with a Lie algebroid structure
such that the anchor is the projection from $A$ to $TM$. A Dirac
structure $L\subset\omni$ is called reducible if $\pomnib(L)$ is a
regular subbundle of $\huaT$. We will see that any Dirac structure
is reducible if $\mathrm{rank}(E)\geq 2$ (Lemma
\ref{Lem:dimention}). The main result is Theorem \ref{Thm:Main},
which claims a one-to-one correspondence between reducible Dirac
structures in $\omni$ and projective Lie algebroids in $\huaT$. In
fact, the projection
 of a reducible Dirac structure $L$ to $\huaT$ yields a projective Lie algebroid $\pomnib(
 L)$ and, conversely,
a projective Lie algebroid $A\subset \huaT$ can be uniquely lifted
to  a Dirac structure $L^A$ using a connection in
 $E$.

 Furthermore,  using the falling operator $\falling{(\cdot)}$,
 we   establish a connection between the derivation algebra $\Der(A)$ of a projective
Lie algebroid $A$   and the normalizer $N_{L^A}$ of the
corresponding lifted Dirac structure $L^A$. We prove that, for any
$X\in N_{L^A}$, $\falling{X}\in \Der(A)$. Conversely, any $\delta
\in \Der(A)$ can be lifted to an element in $N_{L^A}$. Another
observation is that, to any Dirac structure $L\subset \omni$, there
associates a representation of $L$ on $E$, namely
$\rho_L:~L\longrightarrow \dev E$ (Proposition
\ref{pro:representation}). So there is an associated cohomology
group $\mathrm{H}^\bullet(L,\rho_L)$. We will see that the
normalizer of $L$ is related with $\mathrm{H}^1(L,\rho_L)$ and the
deformation of  $L$ is related with $\mathrm{H}^2(L,\rho_L)$.

This paper is  organized as follows. In Section \ref{Sec:omnireview}
we recall  the basic properties of omni-Lie algebroids. In Section
\ref{DiracReduct}, we state  the main result of this paper --- the
correspondence between  reducible Dirac structures and projective
Lie algebroids. In Section 4, several interesting examples are
discussed. In Section \ref{sec:1cohomology}, we study the relation
between the normalizer of a reducible Dirac structure   and  Lie
derivations. In Section \ref{sec:2cohomology}, we give some
applications of the related cohomologies of Dirac structures.

\section{Omni-Lie Algebroids}\label{Sec:omnireview}

We use  the following convention throughout the paper: $E\lon M$
denotes a vector bundle $E$ over a smooth manifold $M$ (we assume
that $E$ is not a zero bundle),  $\dM:~\Omega^\bullet(M)\lon
\Omega^{\bullet+1}(M)$ the usual deRham differential of forms and
$m$
 an arbitrary point in $M$. By $\mathcal{T}$ we denote the direct sum
${~TM\oplus E}$ and use  $pr_{TM}$, $pr_E$, respectively, to denote
the projection from $\mathcal{T}$ to   $TM$ and $E$.

First, we briefly review  the notion of   omni-Lie algebroids
defined in \cite{CLomni}, which generalizes  omni-Lie algebras
defined by A. Weinstein in \cite{Weinomni}. Given a vector bundle
$E$, let ${\jet} E$ be the first jet  bundle of $E$ \cite{Jet
bundle}, and $\dev E$ the gauge Lie algebroid of $E$
\cite{Mkz:GTGA}. These two vector bundles associate, respectively,
with the jet sequence:
\begin{equation}\label{Seq:JetE}
\xymatrix@C=0.5cm{0 \ar[r] & \Hom(TM,E)  \ar[rr]^{\quad\quad\quad\e}
&&
                {\jet}{E} \ar[rr]^{\p} && E \ar[r]  & 0,
                }
\end{equation}
and the Atiyah sequence:
\begin{equation}\label{Seq:DE}
\xymatrix@C=0.5cm{0 \ar[r] & \gl(E)  \ar[rr]^{\id} &&
                \dev{E}  \ar[rr]^{\jd} && TM \ar[r]  & 0.
                }
\end{equation}
The embedding maps $\e$ and $\id$ in the above two exact sequences
will be ignored when there is no risk of confusion. It is well known
that  $\dev{E}$ is a transitive Lie algebroid over $M$, with the
anchor $\jd$ as above \cite{KSK:2002}. The $E$-duality between two
vector bundles is defined as follows.
\begin{defi}Let $A$, $B$ and $E$ be  vector bundles over $M$. We say that
 $B$ is an $E$-dual bundle of
$A$ if  there is a $\CWM$-bilinear   $E$-valued pairing  $
\conpairing{\cdot,\cdot}_E: ~~ A\times_M B\lon E $ which is
nondegenerate, that is, the map $a\mapsto \conpairing{a,\cdot}_E $
is an embedding of $A$ into $\Hom(B,E)$, and similarly for the
$B$-entry.
\end{defi}
  An important result  in
\cite{CLomni} is
 that   $\jet{E}$   is an $E$-dual bundle of
$\dev{E}$ with some nice properties. In fact, we have a
nondegenerate $E$-pairing $\conpairing{\cdot, \cdot}_E$ between
$\jet E$ and $\dev{E}$:
\begin{eqnarray}\nonumber
\conpairing{\mu,\frkd}_E=\conpairing{\frkd,\mu}_E &\defbe& \frkd
u,\quad\forall ~~ \mu=[u]_m\in {\jet
E},~u\in\Gamma(E),~\frkd\in\dev{E}.
\end{eqnarray}
Moreover, this pairing  is $C^\infty(M) $-linear and satisfies the
following properties:
\begin{eqnarray*}
\conpairing{\mu, \Phi }_E &=& \Phi\circ \p(\mu),\quad\forall ~ ~\Phi\in \gl(E),~\mu\in{\jet E};\\
\conpairing{ {\frky} ,\frkd}_{E} &=& {\frky}\circ
\jd(\frkd),\quad\forall ~~ \frky\in \Hom(TM,E),~\frkd\in\dev{E}.
\end{eqnarray*}
An equivalent expression is that we can define $\jet E$ by $\dev E$,
\begin{eqnarray*}  {\jet E}  &\cong &
\set{\nu\in \Hom( \dev{E}  ,E )\,|\,  \nu(\Phi)=\Phi\circ
\nu(\Id_E),\quad\forall ~~ \Phi\in \gl(E )}  \subset \Hom( \dev{E}
,E ).
\end{eqnarray*}
Conversely, $\dev E$ is also determined by $\jet E$:
\begin{eqnarray*}
 {\dev E}   &\cong & \{\delta\in \Hom( {\jet E} ,E )\,|\, \exists~
x\in T M, \ \mbox{ s.t. } \delta(\frky)=\frky(x),\quad\forall ~
~\frky\in \Hom(TM,E) \}.
\end{eqnarray*}

For a Lie algebroid  $(\huaA,[\cdot, \cdot],\alpha)$ over $M$,  a
\emph{representation} of $\huaA$ on a vector bundle $E\lon M$ is a
Lie algebroid morphism $\Rep:~\huaA\lon \dev E$. We may also refer
to $E$ as an $\huaA$-module. To such a representation, there
associates a cochain complex $\sum_{i\geq
0}\Omega^i(\huaA,E)=\sum_{i\geq 0} \Gamma(\Hom(\wedge^i\huaA,E))$
with the coboundary operator:
$$
\dM_{\huaA}: \Omega^\bullet(\huaA, E)\lon \Omega^{\bullet+1}(\huaA,
E),
$$
 defined in a similar fashion as that of the deRham differential \cite{Mkz:GTGA}.
Since $\dev{E}$ is a Lie algebroid and $E $ is a natural
$\dev{E}$-module, we have the cochain complex:
$$\Omega^\bullet(\dev E, E)=\Gamma(
\Hom(\wedge^\bullet{\dev{E}},E))$$  with the coboundary operator:
\begin{equation}\label{Eqt:jetdDEE}
\jetd: \Omega^\bullet(\dev E, E)\lon \Omega^{\bullet+1}(\dev E, E).
\end{equation}
Note that, $\forall$ $u\in\Gamma (E)$, $\jetd u\in \Omega^1(\dev E,
E)$
  is a section of $ \jet E$ and we have a formula:
$$
\jetd(fu)=f\jetd u+ \dM f\otimes u,\quad\forall~ f\in \CWM,~ u\in
\Gamma(E).
$$
The section space $\Gamma (\jet E)$ is an invariant subspace of the
Lie derivative $\LieDerivation_{\frkd}$ for any
  $\frkd \in\Gamma(\dev{E})$. Here $\LieDerivation_{\frkd}$ is
   defined by the
 Leibniz rule as follows:
\begin{eqnarray}\nonumber
\conpairing{\LieDerivation_{\frkd}\mu,\frkd\pie}_{E}&\defbe&
\frkd\conpairing{\mu,\frkd\pie}_{E}-\conpairing{\mu,[\frkd,\frkd\pie]_{\dev}}_{E},
\quad\forall~ \mu \in \Gamma(\jet{E}), ~
~\frkd\pie\in\Gamma(\dev{E}).
\end{eqnarray}

\begin{defi}{\rm\cite{CLomni}}
We call the quadruple
$(\omni,\Dorfman{\cdot,\cdot},\ppairingE{\cdot,\cdot},\omnirho)$ an
{\em omni-Lie algebroid}, where $ \omni=\dev{E}\oplus \jet{E}$,
$\omnirho$  is the projection from $\omni$ to $\dev{E}$, the bracket
$
\Dorfman{\cdot,\cdot}:\Gamma(\omni)\times\Gamma(\omni)\longrightarrow\Gamma(\omni)$
is defined by $$ \Dorfman{\frkd+\mu,\frkr+\nu}\defbe
[\frkd,\frkr]_{\dev}+\LieDerivation_{\frkd}\nu-\LieDerivation_{\frkr}\mu
+ \jetd\conpairing{\mu,\frkr}_E\,, $$ and $\ppairingE{\cdot,\cdot}$
is a nondegenerate symmetric $E$-valued 2-form on $\omni$ defined
by:
\begin{eqnarray*}
\ppairingE{\frkd+\mu,\frkr+\nu}&\defbe&
\half(\conpairing{\frkd,\nu}_E +\conpairing{\frkr,\mu}_E),
\end{eqnarray*}
for any $\frkd,~\frkr\in\dev{E},~\mu,~\nu\in\jet{E}.$
\end{defi}

\begin{thm}\label{Thm:Property of omni}{\rm\cite{CLomni}}
   An omni-Lie algebroid satisfies the following properties:
\begin{itemize}
\item[1)] $(\Gamma(\omni),\Dorfman{\cdot,\cdot})$ is a Leibniz
algebra, \item[2)]
$\omnirho\Dorfman{X,Y}=[\omnirho(X),\omnirho(Y)]_{\dev}$,
 \item[3)]
$\Dorfman{X,fY}=f\Dorfman{X,Y}+(\jd\circ \omnirho(X))(f)Y$,
 \item[4)] $\Dorfman{X,X}= \jetd \ppairingE{X,X}$, \item[5)]
$\omnirho(X)\ppairingE{Y,Z}=\ppairingE{\Dorfman{X,Y},Z}+\ppairingE{Y,\Dorfman{X,Z}}$,
\end{itemize}
for any $X,~Y,~Z\in \Gamma(\omni)$ and $f\in\CWM$.
\end{thm}
It is easy to obtain the following equalities:
\begin{eqnarray}\label{Eqt:DorfmanfXY}
\LDorfman{fX, Y}&=&f\LDorfman{X,Y}-(\jd\circ \omnirho(Y))(f)Y+ 2 \dM
f\otimes \ppairingE{X,Y}, \\
\label{Eqt:DorfmanXYYX} \LDorfman{ X, Y}+\LDorfman{Y,X}&=& 2\jetd
\ppairingE{X,Y}.
\end{eqnarray}
For a subbundle $S\subset \omni$,  we denote
$$
S^\bot= \set{X\in \omni ~|~ \ppairingE{X,s}=0,\quad\forall~ s\in S}.
$$
 We call $S$  isotropic with respect to $\ppairingE{\cdot,\cdot}$ if
$S\subset S^\bot$.

\begin{defi}{\rm\cite{CLomni}}
A Dirac structure in the omni-Lie algebroid $\omni$ is a maximal
isotropic\footnote{One may prove that $L$ is maximal isotropic if
and only if $L=L^\bot$.} subbundle $L\subset \omni$
  such that   $\Dorfman{\Gamma(L),\Gamma(L)}\subset \Gamma(L)$.
\end{defi}

\begin{pro}\label{pro:representation}{\rm\cite{CLomni}}
A Dirac structure $L$ is necessarily a Lie algebroid with the
restricted bracket and the  anchor  $\jd\circ \omnirho$. Moreover,
$\omnirho_L=\omnirho|_L: L \rightarrow \dev{E}$ is a representation
of $L$ on $E$.
\end{pro}
For $\mathcal{T}= {~TM\oplus E}$, we have the standard decomposition
$$\Hom(\mathcal{T},E) = \gl(E)\oplus\Hom(TM,E).$$
 The following exact sequence will be referred as the omni-sequence
of    $E$.
\begin{equation}\label{omni-sequence}
\xymatrix@C=0.5cm{0 \ar[r] & \Hom(\mathcal{T}, E)
\ar[rr]^{\quad\quad\eomnib} &&
                {\omni}  \ar[rr]^{\pomnib} && \mathcal{T} \ar[r]  & 0,
                }
\end{equation}
where the maps $\eomnib$ and $\pomnib$ are defined, respectively, by
\begin{eqnarray*}
\eomnib(\Phi+{\frky})&=&\id(\Phi)+\e({\frky}), \quad\forall ~~
\Phi\in \gl(E), ~{\frky}\in\Hom(TM,E),\\
\pomnib(\frkd+\mu)&=& \jd(\frkd)+\p(\mu),\quad\forall ~~
\frkd\in\dev{E},~\mu\in\jet E.
\end{eqnarray*}
We   regard $\Hom(\mathcal{T} ,E)$ as a subbundle of $\omni$ and
omit the embedding $\eomnib$. Evidently,   $\Hom(\mathcal{T},E)$ is
a maximal isotropic subbundle of $ \omni$. In fact, it is a Dirac
structure of $\omni$ and the bracket is given by
$$
\Dorfman{\alpha,\beta}=\alpha\circ\beta-\beta\circ\alpha,\quad
\forall~\alpha,\beta\in\Gamma(\Hom(\mathcal{T},E)).
$$
In particular, if $\alpha=\Phi+\phi,~\beta=\Psi+\psi$, where
$\Phi,\Psi\in\Gamma(\gl(E)),~\phi,\psi\in\Gamma(\Hom(TM,E))$, then
$$
\Dorfman{\Phi,\Psi}=\Phi\circ\Psi-\Psi\circ\Phi,\quad\Dorfman{\phi,\psi}=0,\quad\Dorfman{\Phi,\phi}=\Phi\circ\phi.
$$
\begin{lem}
\label{lem: idea hom(TM,E)}\begin{itemize}\item[\rm(1)]The subspace
$\Gamma(\Hom(\mathcal{T}, E))$ is a right ideal of $\Gamma(\omni)$.
\item[\rm(2)]For any $h\in \Gamma(\Hom(\mathcal{T},E))$, $X\in
\Gamma(\omni)$, we have
\begin{equation}\label{Eqt:DorfmanhX}
\pomnib\Dorfman{h,X}=h(\pomnib(X)).
\end{equation}\end{itemize}
\end{lem}

\pf For any $  ~X=\frkd+\mu\in \Gamma(\omni)$ and $h=\Phi+{\frky}\in
\Gamma(\Hom(\mathcal{T}, E))$, we have
\begin{equation}\nonumber
\Dorfman{\frkd+\mu,\Phi+{\frky}}=[\frkd,\Phi]_\dev+\LieDerivation_\frkd
\eta-\LieDerivation_\Phi\mu+\jetd\langle\mu,\Phi\rangle_E.
\end{equation}
Since $$\p
(-\LieDerivation_\Phi\mu+\jetd\langle\mu,\Phi\rangle_E)=-\Phi\p(\mu)+\langle\mu,\Phi\rangle_E=0$$
and $\jd[\frkd,\Phi]_\dev$ = $0$, we  have
$$\Dorfman{\frkd+\mu,\Phi+{\frky}} \in\Gamma(\Hom(\mathcal{T},
E)) ,$$ which implies that $\Gamma(\Hom(\mathcal{T}, E))$ is a right
ideal of $\Gamma(\omni)$.

On the other hand,
 we have
\begin{eqnarray*}
\pomnib\Dorfman{h,X}&=&\pomnib([\Phi,\frkd]_{\dev}+\LieDerivation_{\Phi}\mu-\LieDerivation_{\frkd}\frky+\jetd
\langle\frkd,\eta\rangle_E)\\
&=&\Phi(\p\mu)+\frky(\jd\frkd)   = h(\pomnib(X)),
\end{eqnarray*}
which completes the proof. \hfill$\square$

\section{Dirac Structures and Their Reductions}
\label{DiracReduct} Let us first study  some basic properties of
maximal isotropic subbundles of $\omni$. For   any  subbundle
$Q\subset\mathcal{T} $, define:
\begin{eqnarray*}
Q^0&\triangleq&\set{h\in \Hom(\mathcal{T},E)| h(Q)=0 }.
\end{eqnarray*}

\begin{lem}\label{Lem:dimention}If $\mathrm{rank}(E)=r$, $\mathrm{dim}(M)=d$,
then for any maximal isotropic subbundle $L\subset \omni$, we have
\begin{equation}\label{Eqt:dimension}
\mathrm{rank}(L_m)=(1-r)\mathrm{rank}(\pomnib(L_m))+ r(d+r),\quad
\forall~ m\in M.
\end{equation}
Consequently, if $r\geq 2$, both $\pomnib(L)$ and $\pomnib(L)^0$ are
regular   subbundles of, respectively, $\mathcal{T}$ and $ \omni$.
If $r=1$, that is, $E$ is a line bundle, then
$\mathrm{rank}(L)=d+1$.
\end{lem}
\pf Since $L$ is maximal isotropic, or equivalently, $L=L^\bot$, it
is not hard to establish the following exact sequence:
\begin{equation}\label{Seq:maximalhorizontal}
\xymatrix@C=0.7cm{ 0 \ar[r] &  (\pomnib(L_m))^0 \ar[rr]^{{\eomnib }}
&&
                  {L_m} \ar[rr]^{{\pomnib  }} &&   \pomnib(L_m)  \ar[r]  &
                0.}
\end{equation}
Therefore, we have
\begin{eqnarray*}
\mathrm{rank}(L_m)&=&\mathrm{rank}(\pomnib(L_m))+\mathrm{rank}(\pomnib(L_m))^0\\
&=&\mathrm{rank}(\pomnib(L_m))+(r+d-\mathrm{rank}(\pomnib(L_m)))\times
r\\
&=&(1-r)\mathrm{rank}(\pomnib(L_m))+ r(d+r).
\end{eqnarray*}
\hfill$\square$
\begin{defi}
For a vector subbundle $A\subset \mathcal{T}$, a section
$s:~A\longrightarrow\omni$ {\rm(i.e. $\pomnib\circ s= \Id_A$)} is
called
 isotropic if its image $s(A)\subset\omni$ is isotropic. Two isotropic sections $s_1$
 and $s_2$ are said to be equivalent
 if $(s_1-s_2)(A)\subset A^0$. The
equivalence class of an isotropic section $s$ is denoted by
$\widetilde{s}$.
\end{defi}

\begin{pro}\label{thm:characteristic pair}
If {\em rank}$ E\geq 2$, there is a one-to-one correspondence
between  maximal isotropic subbundles $L \subset\omni$ and pairs
$(A,\widetilde{s})$, where $A$ is a subbundle of $\huaT$ and
$s:~A\lon \omni$ is an isotropic section.
\end{pro}

For this reason, we call   $(A,\widetilde{s})$ the {\em
characteristic pair} of
  $L$, and  write $L=L_{s,A}$.

\pf Let  $L\subset\omni$ be a maximal isotropic subbundle and
$A=\pomnib(L)$. By Lemma \ref{Lem:dimention}, $A$ is a regular
subbundle. Any split $s:~A\lon L$ of the corresponding exact
sequence (\ref{Seq:maximalhorizontal}) yields an isotropic section
and $(A,\widetilde{s})$ is defined to be the characteristic pair of
$L$. It is well defined  since for any two isotopic sections $s_1$,
$ s_2$, we have $\Img(s_1-s_2)\subset \pomnib(L)^0=A^0$, which is
equivalent to  $\widetilde{s_1}=\widetilde{s_2}$.

Conversely,  given a subbundle $A\subset \mathcal{T}$ and any
characteristic pair $(A,\widetilde{s})$, set $L_{s,A}= s(A)\oplus
A^0$. Evidently, $L_{s,A}$ is a maximal isotropic subbundle of
$\omni$ whose characteristic pair is $(A,\widetilde{s})$.  It is
also clear that
  if $\widetilde{s_1}=\widetilde{s_2}$,  $L_{s_1,A}=L_{s_2,A}$.

 One may check that these two constructions
are inverse to each other.    \hfill$\square$

\begin{defi}\label{Defi:projective Lie algebroid}
A  projective Lie algebroid  is a subbundle  $ A \subset TM\oplus E$
which is a Lie algebroid  $(A,
 [ \cdot, \cdot ]_A, \rho_A)$ and the anchor  $\rho_A = pr_{TM}|_A$.
 \end{defi}

\begin{ex}\rm   Let
 $\huaA\longrightarrow N$  be a Lie algebroid over a smooth manifold $N$
 and $\alpha$ its anchor. Let
$f:M \longrightarrow N$ be a smooth map and  $f^*\huaA \rightarrow
M$  the pull back bundle along $f$. We denote the pull back Lie
algebroid of $\huaA$ over $M $ by $f^{! }A=TM \oplus_{TN}\huaA$,
which is given by
$$
TM \oplus_{TN}\huaA=\set{(x ,X)\in T_{m }M \oplus \huaA_{f(m )}| m
\in M ,\mbox{ and } f_*(x )=\alpha(X)}.
$$
 Sections of
$TM \oplus_{TN}\huaA$ are of the form:
$$ x \oplus(\sum u_i \otimes
X_i),\qquad x \in\frkX(M ),~u_i \in C^\infty(M ),~X_i\in
\Gamma(\huaA),
$$
such that $f_*(x  (m ))=\sum u_i (m )\alpha(X_i(f(m )))$. The anchor
$\alpha^!$ of the Lie algebroid $f^{! }\huaA$ is   the projection to
the first summand. The Lie bracket can be \emph{locally} expressed
by
\begin{eqnarray*}
&&[x \oplus(\sum u_i \otimes X_i),y \oplus(\sum
v_j \otimes Y_j)]\\
&=&[x ,y ]\oplus(\sum u_i  v_j \otimes [X_i,Y_j]+\sum x (v_j
)\otimes Y_j-\sum y (u_i )\otimes X_i).
\end{eqnarray*}
Thus   the pull back Lie algebroid $f^{! }\huaA$ of the Lie
algebroid $\huaA$ is a projective Lie algebroid in $TM\oplus
f^*\huaA$.
\end{ex}
\begin{ex}\rm We suppose that the base manifold $M$ is
compact and let $H\subset TM$ be an integrable distribution. It is
well known that there is some vector  bundle $E$ such that the
vector bundle $F=H\oplus E$ is trivial. Suppose that $\rank F=n$ and
$\varepsilon_1,\cdots,\varepsilon_n$
 are everywhere linear independent sections of $F$, i.e. a frame of
$\Gamma(F)$. Write $\varepsilon_i=x_i+e_i$, where $x_i$ and $e_i$
are sections of $H$ and $E$ respectively. It is clear that
$\Gamma(H)=$span$\{x_1,\cdots,x_n\}$ and
$\Gamma(E)=$span$\{e_1,\cdots,e_n\}$ (over $\CWM$). Since $H$ is an
integrable distribution, there exist functions $c_{i,j}^k\in
C^\infty(M)$ such that $[x_i,x_j]=c_{i,j}^kx_k$. Now set
$[\varepsilon_i,\varepsilon_j]=c_{i,j}^k\varepsilon_k$. It is easy
to see that $F$ is a projective Lie algebroid in $TM\oplus E$.
\end{ex}

A Dirac structure $L\subset\omni$ is called $\mathbf{reducible}$ if
$\pomnib(L)$ is a regular subbundle of $\huaT$. By Lemma
\ref{Lem:dimention}, any Dirac structure is reducible if
$\mathrm{rank}(E)\geq 2$. As a main result of this paper, the
following theorem   describes the nature of reducible Dirac
structures in the omni-Lie algebroid $\omni$.
\begin{thm}\label{Thm:Main}
For any vector bundle $E$,  there is a one-to-one correspondence
between reducible Dirac structures $L \subset \omni$  and projective
Lie algebroids $A=\pomnib(L) \subset \huaT$ such that $A$ is the
quotient Lie algebroid of $L$.
\end{thm}

\pf Assume that $L$ is a reducible Dirac structure and let
$A=\pomnib(L)\subset \mathcal{T}$. Then we have the following exact
sequence:
\begin{equation}\label{Dirac-sequence2}
\xymatrix@C=0.7cm{ 0 \ar[r] & A^0 \ar[rr]^{{\eomnib }} &&
                  {L} \ar[rr]^{{\pomnib }} && A \ar[r]  &
                0.}
\end{equation}
By $L$ being reducible,  $A$ is a  regular subbundle,  $A^0$ as
well. The anchor $\jd\circ \omnirho$ vanishes if restricted on
$A^0$. Furthermore,
 by Lemma \ref{lem: idea hom(TM,E)} and the fact that $L$ is a
Dirac structure,   $A^0$ is an ideal of $L$. So we have a quotient
Lie algebroid structure $(A,[\cdot,\cdot]_A,\rho_A)$, where $\rho_A$
is clearly the projection to $TM$. This proves that $A$ is indeed a
projective Lie algebroid.

Conversely,  for the projective Lie algebroid
$(A,[\cdot,\cdot]_A,\rho_A)$, define a subset $L^A\subset
\pomnib\inverse(A)\subset \omni$ by:
\begin{eqnarray}\label{L}\nonumber
L^A_m&\defbe& \{X\in \pomnib\inverse(A)_m|\mbox{ for some
}\widetilde{X}\in\Gamma (\pomnib\inverse(A))\mbox{ with
}\widetilde{X}_m=X, \quad\mbox{there holds}\\\label{Set:L}
&&\quad\quad\quad\quad
\pomnib\Dorfman{\widetilde{X},Y}_m=([\pomnib\widetilde{X},\pomnib
Y]_{A})_m\,,\quad\forall ~~ Y\in\Gamma( \pomnib\inverse(A))\}.
\end{eqnarray}
Note that by Equation (\ref{Eqt:DorfmanfXY}), we have
$$\pomnib\Dorfman{f\widetilde{X},Y}_m-([f\pomnib\widetilde{X},\pomnib
Y]_{A
})_m=f(\pomnib\Dorfman{\widetilde{X},Y}_m-([\pomnib\widetilde{X},\pomnib
Y]_{A })_m).
$$
Hence the above definition does not depend on the choice of
$\widetilde{X}$.

To prove that $L^A$ is the unique reducible Dirac structure such
that the induced projective Lie algebroid is
$(A,[\cdot,\cdot]_A,\rho_A)$, we need three steps as follows.  Step
1
 proves that $L^A$ is a maximal isotropic subbundle such that
$\pomnib(L^A)=A$.  Step 2 proves that $L^A$ is closed under the
bracket $\Dorfman{\cdot,\cdot}$ and it follows that $L^A$ is a
reducible Dirac structure such that the induced projective Lie
algebroid is $(A,[\cdot,\cdot]_A,\rho_A)$. The last step proves the
uniqueness of such Dirac structures.

{\bf Step 1.} We prove that $L^A$ is a maximal  isotropic subbundle.
We will construct a maximal isotropic subbundle $L_{s_\gamma,A}$
 using a connection  $\gamma $ in the vector bundle $E$ and prove that $L_{s_\gamma,A}=L^A$.

 Recall that a connection in $E$ is a bundle map
 ${\connection}:~TM~\rightarrow~\dev E$ such that
 $\jd\circ\connection=\Id_{TM}$. Associated with $\connection$ there is a
 back connection $\omega:~ \dev E\lon\gl(E)$, such that
  $\id\circ\omega+{\connection}\circ\jd=\Id_{\dev E}$.
So we can define a bundle map $\widetilde{{\connection}}:~ E\lon
\jet{E}$ by
\begin{equation}\label{defi:gamma tuta}
\conpairing{\widetilde{{\connection}}(e),\frkd}_E\defbe
\omega(\frkd)(e) =(\frkd-\connection\circ\jd(\frkd))(e),\quad\forall
~ \frkd\in \dev E
\end{equation} such that $\p\circ\widetilde{{\connection}}=\Id_{E}$.
In turn, we get a map:
\begin{equation}\label{eqt:gamma}
{\connection}+\widetilde{{\connection}}:~ \mathcal{T}\lon \omni
\quad\mbox{ such that}\quad
\pomnib\circ({\connection}+\widetilde{{\connection}})=\Id_{\mathcal{T}}.
\end{equation}
We still denote   this map by ${\connection}$. This does not make
any confusion since it depends on what is put right after it.

Choose an arbitrary subbundle $C\subset\mathcal{T}$, such that
$\mathcal{T}=A\oplus C$. Define a bundle map $\Omega_{\connection}:
~\mathcal{T}\wedge \mathcal{T}\lon E$ by
\begin{eqnarray*}
\Omega_{\connection}(a,b)&=&
[a,b]_A-\pomnib\Dorfman{{\connection}(a),{\connection}(b)},\quad
\forall ~~ a,b\in \Gamma(A),\\
\Omega_{\connection}(c,t)&=&0, \qquad \forall ~c\in C,~t\in
\mathcal{T}.
\end{eqnarray*}
To see that $\Omega_{\connection}\in
 \Hom(\wedge^2\mathcal{T},E)$, first for any
$a=x+u,~b=y+v\in \Gamma(A)$, where $x,y\in \frkX(M)$,
$u,v\in\Gamma(E)$, we have
\begin{eqnarray*}
 \pomnib\Dorfman{{\connection}(x+u),{\connection}(y+v)}
&=&\pomnib([{\connection}(x),{\connection}(y)]_{\dev}+\LieDerivation_{{\connection}(x)}{\connection}(v)-
\LieDerivation_{{\connection}(y)} {\connection}(u)+\jetd\langle{\connection}(y),{\connection}(u)\rangle_E)\\
&=&[\jd{\connection}(x),\jd{\connection}(x)]_{\dev}+{\connection}(x)(\p{\connection}(v))-
{\connection}(y)(\p{\connection}(u))
\\ &=&[x,y]+{\connection}(x)v-{\connection}(y)u,
\end{eqnarray*}
which implies that
\begin{equation}\label{Eqt:Omegagamma}
\Omega_{\connection}(x+u,y+v)=([x+u,y+v]_A-[x,y])-{\connection}(x)(v)+{\connection}(y)(u).
\end{equation}
Thus we have $\Omega_{\connection}(x+u,y+v)\in\Gamma(E)$. On the
other hand, for any $f\in C^\infty
 (M)$, we have
\begin{eqnarray}\label{11}
\nonumber\Omega_{\connection}(x+u,f(y+v))&=&([x+u,f(y+v)]_A-[x,fy])-{\connection}(x)(fv)+{\connection}(fy)(u)\\
&=&\nonumber f\Omega_{\connection}(x+u,y+v)+x(f)(y+v)-x(f)y-\alpha({\connection}(x))(f)v\\
&=&\label{eqt:omega f} f\Omega_{\connection}(x+u,y+v).
\end{eqnarray}
By (\ref{Eqt:Omegagamma}) and (\ref{11}), we obtain that
$\Omega_{\connection}\in
 \Hom(\wedge^2\mathcal{T},E)$.
 We also denote the associated skew-symmetric
 map from $\mathcal{T}$ to $\Hom(\mathcal{T},E)$ by
 $\Omega_{\connection}$.

Define an isotropic section $s_\gamma:~A\longrightarrow \omni$ by
$$s_\gamma(a)=\gamma(a)+\Omega_\gamma(a),\quad\forall ~ a\in A.$$
 In fact, for $a=x+u,~b=y+v\in
\Gamma(A)$, we have
\begin{eqnarray*}
&&\ppairingE{s_{\connection}(x+u),s_{\connection}(y+v)}\\&=&
\ppairingE{\connection(x)+\connection{(u)}+\Omega_{\connection}(a),\connection(y)+\connection{(v)}+\Omega_{\connection}(b)}
\\
&=&\half(\Omega_{\connection}(y+v,x+u)+\Omega_{\connection}(x+u,y+v)
+\conpairing{\gamma(x),\gamma(v)}_E+\conpairing{\gamma(y),\gamma(u)}_E
)=0.
\end{eqnarray*}
By Proposition \ref{thm:characteristic pair}, we get a maximal
isotropic subbundle $L_{s_\gamma,A}$:
\begin{equation}\label{eqt:construction Dirac L}
L_{s_\gamma,A}=\gamma(A)+\Omega_\gamma(A)+A^0.
\end{equation}
We can directly check that $L_{s_\gamma,A}$ does not depend on the
choice of the connection $\gamma$ and the subbundle $C$. An
alternate approach is to  prove that  $L_{s_\gamma,A}=L^A$, since
$L^A$ does not depend on $s_\gamma$ and $A$.

Now we prove $L_{s_\gamma,A}=L^A$. Any $X\in \Gamma(L_{s_\gamma,A})$
 has the form
 $X={\connection}(a)+\Omega_\gamma(a)+h,$
where $a=x+u\in\Gamma(A)$ and $h\in\Gamma(A^0)$. For any
$Y=\frkd+\mu\in\Gamma(\pomnib\inverse(A))$ satisfying
$\pomnib(Y)=y+v\in\Gamma(A)$, we have
\begin{eqnarray*}
\pomnib\Dorfman{X,Y}&=&
\pomnib(\Dorfman{\connection(x)+\connection(u),\frkd+\mu}+\Dorfman{\Omega_\connection(a)+h,Y})\\
&=&
\pomnib([{\connection}(x),\frkd]_{\dev}+\LieDerivation_{{\connection}(x)}\mu
-\LieDerivation_{\frkd}{\connection}(u)+\jetd\langle\gamma(u),\frkd\rangle_E )+(\Omega_\connection(a)+h)(\pomnib(Y))\\
&=&[x,\jd\frkd]
+{\connection}(x)(v)-\frkd(u)+\langle\gamma(u),\frkd\rangle_E
+\Omega_{{\connection}}(x+u,y+v)\\
&=&[x,y]+{\connection}(x)v-{\connection}(y)u+\Omega_{{\connection}}(x+u,y+v)\\&=&[x+u,y+v]_A,
\quad\mbox{(using (\ref{Eqt:Omegagamma}) ) }\\
&=&[\pomnib(X),\pomnib(Y)]_A.
\end{eqnarray*}
Thus, $X\in \Gamma(L^A)$. So we have $L_{s_\gamma,A}\subset L^A$.
 Since $\pomnib(L^A)\subset A$, any
 $X\in L^A$ can be written as $X=X_0+h$, where
$X_0\in L_{s_\gamma,A}$ and $h\in \Hom(\huaT,E)$. Thus $h=X-X_0\in
L^A\cap \Hom(\huaT,E)$.

For any $ k $ $\in \Hom(\huaT_m,E_m)$ = $\Ker \pomnib$ and  $
\widetilde{k}\in \Gamma(\Hom(\huaT ,E) )$ satisfying
$\widetilde{k}(m)=k$, $\forall$ $ Y \in\Gamma (\pomnib\inverse(A))$,
we have, by Equation (\ref{Eqt:DorfmanhX})
\begin{eqnarray*}
\pomnib\Dorfman{ \widetilde{k} ,Y}_m-([\pomnib(\widetilde{k}
),\pomnib(Y)]_A)_m &=& k(\pomnib(Y)).
\end{eqnarray*}
Thus $k\in L^A_m\cap\Hom(\huaT_m,E_m)$ if and only if $k\in A_m^0$,
that is,
\begin{equation}
\label{relationtemp1}L^A\cap \Hom(\huaT,E)= A^0.
\end{equation}
So we have proved that $L^A\subset L_{s_\gamma,A}$. By maximality,
$L^A=L_{s_\gamma,A}$ and hence $L^A$ is a maximal isotropic
subbundle of $\omni$.

{\bf Step 2.} We prove that $\Gamma(L^A)$ is closed under the
 bracket operation $\{\cdot,\cdot\} $ and it follows that $L^A=L_{s_\gamma,A}$
is a reducible Dirac structure.

 For any  $X_1$, $X_2\in\Gamma (L^A)$ and $Y\in
\Gamma (\pomnib\inverse(A))$, we have  $\Dorfman{X_1,X_2}\in \Gamma
(\pomnib\inverse(A))$ and  $\Dorfman{X_i,Y}\in\Gamma
(\pomnib\inverse(A))$. Moreover, we have
\begin{eqnarray*}
\pomnib\Dorfman{\Dorfman{X_1,X_2},Y}
&=&\pomnib\Dorfman{X_1,\Dorfman{X_2,Y}}-
\pomnib\Dorfman{X_2,\Dorfman{X_1,Y}}\\
&=&[\pomnib X_1,\pomnib\Dorfman{X_2,Y}]_A
-[\pomnib X_2,\pomnib\Dorfman{X_1,Y}]_A\\
&=&[\pomnib X_1,[\pomnib X_2,Y]_A]_A -[\pomnib X_2,[\pomnib
X_1,Y]_A]_A\\
&=&[[\pomnib X_1,\pomnib X_2]_A,\pomnib Y]_A\\&=&
[\pomnib\Dorfman{X_1,X_2},\pomnib Y]_A,
\end{eqnarray*}
which implies that $\Dorfman{X_1,X_2}\in \Gamma (L^A)$. So $L^A$ is
a Dirac structure. In Step 1, we have proved that $\pomnib(L^A)=A$,
and in turn, $L^A$ is a reducible Dirac structure. By definition,
the induced projective Lie
  algebroid  is  exactly $(A,[\cdot~,~\cdot]_A,\rho_A)$.

 {\bf Step 3.}  We prove the uniqueness of such Dirac structures.

Assume that   ${L'}$ is  another reducible Dirac structure
satisfying the same requirements. It suffices to  prove that
${L'}\subset L^A$, since $L^A$ is a maximal isotropic subbundle. For
any $X \in {{L'}_m}$ and $\widetilde{X}\in\Gamma ({L'})$ such that
$\widetilde{X}_m=X$, we prove that $X\in L^A_m$. In fact, $\forall$
$Y\in \Gamma (\pomnib\inverse(A))$, we are able to find some $Y'\in
\Gamma ({L'})$ such that $\pomnib Y'=\pomnib Y$. So we can write
$Y=Y'+K$, where $K\in \Gamma (\Hom(\huaT,E))$. By Lemma \ref{lem:
idea hom(TM,E)}, $\Dorfman{\widetilde{X},K}\in \Gamma
(\Hom(\huaT,E))$. Thus,
\begin{eqnarray*}
\pomnib\Dorfman{\widetilde{X},Y}=\pomnib\Dorfman{\widetilde{X},Y'}
+\pomnib\Dorfman{\widetilde{X},K} =[\pomnib\widetilde{X},\pomnib
Y']_A=[\pomnib\widetilde{X},\pomnib Y ]_A,
\end{eqnarray*}
which implies  that $X\in L^A_m$. So we have ${L'}\subset L^A$. The
proof of Theorem \ref{Thm:Main} is thus completed. \hfill$\square$

The projective Lie algebroid $\pomnib(L)$ is called the
$\mathbf{reduction}$ of the reducible Dirac structure $L$. The
reducible Dirac structure $L^A$ is called the $\mathbf{lift}$ of the
projective Lie algebroid $A$.

\section{Some Examples}
Bellow we give some basic examples of Dirac structures in the
omni-Lie algebroid.
\begin{ex}\rm
For a vector space $V$,   our theorem claims a one-to-one
correspondence between  Dirac structures of the omni-Lie algebra $
\gl(V)\oplus V$ and Lie algebra structures on subspaces of $V$. Thus
Dirac structures characterize not only all   Lie algebra structures
on $V$, as pointed out by Weinstein \cite{Weinomni}, but also  all
  Lie algebra structures on
  subspaces of $V$.
 \end{ex}
\begin{ex}\rm\label{Ex:Graphmu}
Given a skew-symmetric bundle map
 ${\widehat{\lambda}}: ~\dev E\lon \jet
E$,  its graph
$$L^{\widehat{\lambda}}=\set{\frkd+{\widehat{\lambda}}(\frkd)~|~\forall ~ \frkd\in\dev E}
\subset \omni
$$
 is clearly a maximal isotropic subbundle. Furthermore, we
have ${\widehat{\lambda}}(\gl(E))\subset\Hom(TM,E)$, i.e.
$\p{\widehat{\lambda}}(\Phi)=0$. In fact, $\forall$ $\Phi\in
\gl(E)$, we have
$$
\conpairing{{\widehat{\lambda}}(\Phi), \Id_E
}_E=\p{\widehat{\lambda}}(\Phi), \qquad
\conpairing{{\widehat{\lambda}}(\Id_E ), \Phi}_E
=\Phi\circ\p{\widehat{\lambda}}(\Id_E ).
$$
Since $\widehat{\lambda}$ is skew-symmetric, we have
$\p{\widehat{\lambda}}(\Phi)=-\Phi\circ\p{\widehat{\lambda}}(\Id_E
)$. If we take $\Phi=\Id_E$, then $\p{\widehat{\lambda}}(\Id_E )=0$.
Thus, $\p{\widehat{\lambda}}(\Phi)=0$.

Let $\lambda:TM\longrightarrow E$ be the induced bundle map of
$\widehat{\lambda}$. Then we have the following commutative diagram:
$$
\xymatrix{
 0 \ar[r] & \gl(E )  \ar[d]_{-\lambda^*} \ar[rr]^{\id} &&
                 \dev{E}  \ar[d]_{{{\widehat{\lambda}}}}  \ar[rr]^{\jd} && T M \ar[d]_{{\lambda}} \ar[r]  & 0\\
 0 \ar[r] & \Hom(T M,E )    \ar[rr]^{\e} &&
                 \jet{E}   \ar[rr]^{\p } && E
                  \ar[r]  & 0.
                }
$$
Here $-\lambda^*$ is induced by ${{\widehat{\lambda}}}|_{\gl(E)}$,
which is given by $\Phi\mapsto -\Phi\circ \lambda$. So we have the
following exact sequence:
$$
\xymatrix@C=0.7cm{ 0 \ar[r] & \Graph_{-\lambda^*} \ar[rr]^{{  }} &&
                  {L^{\widehat{\lambda}}} \ar[rr]^{{  }} && \Graph_{\lambda} \ar[r]  &
                0,}
$$
where $\Graph_{\lambda}=\pomnib(L^{\widehat{\lambda}})$ is the graph
of $\lambda$ and $\Graph_{-\lambda^*}=L^{\widehat{\lambda}}\cap
\Hom(\huaT,E)$ is the graph of ${-\lambda^*}$.
 \end{ex}


We claim that the following three  statements are equivalent.
\begin{itemize}
\item[1)]
$L^{\widehat{\lambda}}$ is a Dirac structure. \item[2)] $\jetd
{\widehat{\lambda}}=0$, regarding $\widehat{\lambda}$ as a map $\dev
E\wedge \dev E\lon E$ in the obvious sense:
$$\widehat{\lambda}(\frkd,\frkr)=\conpairing{\widehat{\lambda}(\frkd),\frkr}_E,\quad
\forall~~ \frkd,\frkr\in \dev E.$$
\item[3)]${\widehat{\lambda}}=-\jetd (\lambda\circ\jd)$.
\end{itemize}

In fact,  1) $\Longleftrightarrow$ 2) is merely some calculations.
 3) $\Longrightarrow$ 2) is trivial. To see the reverse, notice that
 $\forall$
$\frkr,\frks\in \Gamma(\dev E)$,
\begin{eqnarray*}
 \jetd \widehat{\lambda}(\Id_E,\frkr,\frks)
&=&\conpairing{\widehat{\lambda}(\frkr),\frks}_E -\frkr
\conpairing{\widehat{\lambda}(\Id_E),\frks}_E+ \frks
\conpairing{\widehat{\lambda}(\Id_E),\frkr}_E
-\conpairing{\widehat{\lambda}[\frkr,\frks]_\dev,\Id_E}\\
&=& \widehat{\lambda}(\frkr,\frks)
+\frkr(\lambda\circ\alpha(\frks))-\frks(\lambda\circ\alpha(\frkr))
-(\lambda\circ\alpha)[\frkr,\frks]_\dev\,,
\end{eqnarray*}
which implies that 2) $\Longrightarrow$ 3).

  Thus,  any Dirac
  structure of the type $L^{\widehat{\lambda}} $ is a reducible Dirac structure and totally determined by
  $$\pomnib(L^{\widehat{\lambda}})=\Graph_{\lambda}\subset
  \huaT ,$$ which is isomorphic to $TM$ and equipped with the induced Lie algebroid
  structure.

\begin{ex}\rm
\rm\label{Ex:Graphpi}(See \cite{CLomni}) For a skew-symmetric bundle
map $\pisharp: \jet{E}\lon \dev{E}$, its graph
$$
L_{\pisharp}=\set{ \pisharp(\mu)+\mu ~|~\mu\in \jet{E}} \subset
\omni
$$
 is clearly a maximal isotropic subbundle of
$\omni$.  It can be proved that $L_{\pisharp}$ is a Dirac structure
if and only if the following equation holds for all
$\mu,\nu\in\Gamma(\jet{E})$,
\begin{equation}\nonumber
\pisharp\pibracket{\mu,\nu}=[\pisharp(\mu), \pisharp(\nu)]_{\dev},
\end{equation}
where the $\pi$-bracket $\pibracket{\cdot,\cdot}$ on
$\Gamma(\jet{E})$ is given by:
\begin{equation}\label{pibracket}
\pibracket{\mu,\nu}\defbe
\LieDerivation_{\pisharp(\mu)}\nu-\LieDerivation_{\pisharp(\nu)}\mu-
\jetd\circ\pi(\mu\wedge \nu).
\end{equation}

To see what   $\pi$ encodes, we need to consider the following two
situations:

$\bullet$ $\mathrm{rank}(E)\geq 2$. In this case,  in \cite{CLomni},
we proved that  such Dirac structures are in one-to-one
correspondence with Lie algebroid structures on $E$. Let us see how
 Theorem \ref{Thm:Main} recovers this result.
On one hand, there is an obvious  one-to-one correspondence between
Lie algebroid structures $(E,[~\cdot,\cdot~]_E,\rho_E)$   and
projective Lie algebroids
  $\Graph_{\rho_E}$ which are the graphs of  $\rho_E: ~E\lon
TM$. On the other hand, by Lemma \ref{Lem:dimention}, any Dirac
structure is reducible. Especially, for any Dirac structure
$L_\pi\subset \omni$, $\pomnib(L_\pi)$ should be a projective Lie
algebroid. However, $\pomnib(L_\pi)$ is also a graph and hence there
is an induced Lie algebroid structure on $E$. So we conclude that
Lie algebroid structures on $E$ are   in one-to-one correspondence
with
  Dirac structures of the type $L_{\pisharp}$.

$\bullet$  $\mathrm{rank}(E)=1$. For any reducible Dirac structure
$L_\pi\subset\omni$, $\pomnib(L_\pi)$ is a projective Lie algebroid.
But in general, it may not be a graph and so there  is no induced
Lie algebroid structure on $E$. However, there is always a local Lie
algebra structure on $E$ associated with the Dirac structure $L_\pi$
(not necessarily reducible) as proved in \cite{CLomni}.
\end{ex}

\begin{ex}\rm\label{Ex:prjectlinebundle}
Consider the  case that $ A\subset \huaT$ is an arbitrary line
bundle, which is naturally  a projective Lie algebroid. In fact, for
any neighborhood  $\huaU\subset M$ such that $A|_{\huaU}$ is
trivial, i.e. there is a nowhere singular section $a=x+u$,  the Lie
bracket of $\Gamma(A|_{\huaU})$ is given by:
$$
[fa,ga]_A=(f x(g)-gx(f))a,\quad\forall~ f,~g\in C^\infty(\huaU).
$$
It is easy to check that this bracket is well defined.

The lifted Dirac structure declared by Theorem \ref{Thm:Main} can be
constructed by Equation (\ref{eqt:construction Dirac L}). Just take
any connection $\gamma$. Since $A$ is a line bundle, we have
$\Omega_\gamma(A)\subset A^0$. The lifted Dirac structure is given
by $L^A=L_{s_\gamma,A}=\gamma(A)\oplus A^0$.

\end{ex}

\begin{ex}\rm
Assume that $F\subset E$ is a vector subbundle and
 $(F,[\cdot,\cdot],\rho_F)$ is a Lie algebroid. Then $\Graph_{\rho_F}$,
 the graph of $\rho_F$ is a projective Lie algebroid. Now we
 construct the lifted Dirac structure. Evidently, we have
 $$
 \Graph_{\rho_F}^0=\set{ \Phi+\frky\in \Hom(\huaT,E)~|~
 (\frky\circ\rho_F+\Phi)|_{F}=0}.
 $$
Let $L_1\subset \omni$ be the subset  generated by elements of the
form
  $\frkd^v_m+[v]_m$, where $m\in M$, $v\in
 \Gamma(F)$, $ \frkd^v_m\in (\dev E)_m$ and they satisfy
 $$ \frkd^v_m(u)=
 ([v,u]_F)_m ,\quad\forall u\in \Gamma(F).$$
Let $L=L_1 + \Graph_{\rho_F}^0$. Accordingly, we get an
 exact sequence:
 $$\xymatrix@C=0.7cm{ 0 \ar[r] &  \Graph_{\rho_F}^0 \ar[rr]^{ }
&&
                  {L } \ar[rr]^{{\pomnib }} && \Graph_{\rho_F}  \ar[r]  &
                0.}$$
It is clear that $L$ is  maximal isotropic. Moreover, for all
$u,v\in \Gamma(F)$, we have $$[\frkd^u,\frkd^v]_{\dev
E}=\frkd^{[u,v]_F}$$ and it follows that  $\Gamma(L)$ is closed
under the  bracket   $\Dorfman{\cdot,\cdot}$. Hence $L$ is the
lifted Dirac structure.
 \end{ex}

\begin{ex}\label{Example:Full-Dirac}\rm
We consider a projective Lie algebroid $A$ which is
\emph{transitive}, i.e. $\rho_A(A)=pr_{TM}(A)=TM$. One can construct
a map $\vartheta: ~TM\lon E$ such that $A=\Graph_{\vartheta}\oplus
E_0$, where $\Graph_{\vartheta}$ is the graph of $\vartheta$ and
$E_0$ is a subbundle of $E$.   In this case, $E_0$ must be a Lie
algebra bundle and we denote its Lie bracket by $[\cdot,\cdot]^0$.

For any vector field  $x\in \frkX(M)$,  write
$\bar{x}=x+\vartheta(x)\in \Gamma(A)$. There is a suitable
connection ${\connection}: TM\lon \dev E_0$ such that
$$[\bar{x},u]_A={\connection}(x)u,\quad\forall ~ x\in \frkX(M),~ u\in \Gamma(E_0).$$
Define ${R}: \wedge^2 TM\lon E_0$ by
$$
{R}(x,y)=[\bar{x},\bar{y}]_A-\overline{[x,y]},\quad\forall~
x,y\in\frkX(M).
$$
So  the Lie bracket of $\Gamma(A)$ can be written as
$$
[\bar{x}+u,\bar{y}+v]_A=\overline{[x,y]}
+{R}(x,y)+{\connection}(x)v-{\connection}(y)u+[u,v]^0,\quad \forall~
\bar{x}+u ,\bar{y}+v\in\Gamma(A).
$$
Under the structure defined by the given data $({\connection},{R})$,
$ A= \Graph_{\vartheta}\oplus E_0$  is a Lie algebroid if and only
if $\forall$ $
  x,y,z\in \frkX(M),~ u,v\in\Gamma(E_0)$, the following
  compatibility
conditions hold
\begin{eqnarray*}
[{\connection}(x)u,v]^0+[u,{\connection}(x)v]^0 &=&{\connection}(x)[u,v]^0,\\
{[{\connection}(x),{\connection}(y)]_{\dev}}-{\connection}[x,y]&=&
\ad^{E_0}_{{R}(x,y)},\\
{R}([x,y],z)-{\connection}(x){R}(y,z)+~c.p.~&=& 0.
\end{eqnarray*}
We extend the connection $\gamma$ in the vector bundle $E_0$ to a
connection $\widetilde{\gamma}$ in the vector bundle $E$. By
(\ref{Eqt:Omegagamma}), we have
\begin{eqnarray*}
\Omega_{\widetilde{\connection}}(\bar{x}+u,\bar{y}+v)
&=&{R}(x,y)+[u,v]^0-(d^{\widetilde{\gamma}} \vartheta)(x,y)\,,
\end{eqnarray*}
where
$$(d^{\widetilde{\gamma}} \vartheta)(x,y)=\widetilde{\gamma}(x)\vartheta(y)
-\widetilde{\gamma}(y)\vartheta(x)-\vartheta[x,y].
$$
The lifted Dirac structure,  given by Theorem \ref{Thm:Main}, can be
expressed by (\ref{eqt:construction Dirac L}).

In particular, if $A=\huaT=TM\oplus E$ (so  that we may take
$\vartheta=0$), then
$$
L_{s_\connection,\huaT}=\set{ {\connection}(x+u)+i_x{R}+\ad^E_u\mid
\forall~  ~ x+u\in {~TM\oplus E }}.
$$

 We note that the above construction of projective Lie algebroids
 includes a standard
type of Lie algebroids, known as a \emph{semi-direct-sum}.
 If $E$ is a vector bundle over $M$ which admits a \textbf{flat}
connection $\nabla:~TM\lon \dev E$, then the direct sum $TM\oplus E$
has a Lie algebroid structure over $M$, for which the   anchor is
the projection to $TM$ and the Lie bracket is given by:
$$
[x+u,y+v]\defbe [x,y]+\nabla_xv-\nabla_yu, \quad\forall~ x+u,y+v\in
\Gamma(TM\oplus E).
$$
\end{ex}

\section{The Normalizer  of  Dirac Structures}\label{sec:1cohomology}
In this   and the next section, we always assume that  Lie
algebroids under consideration are not zero. For a  Lie algebroid
$A$,  call $\Der(A)$,  the set of Lie derivations of $A$:
$$
\Der(A)= \set{\delta\in \Gamma(\dev A)~|~ \delta[a_1,a_2]_A=[\delta
a_1,a_2]_A+[a_1,\delta a_2]_A,\quad \forall ~ a_1,a_2\in \Gamma(A),
}
$$
the derivation algebra of $A$.
\begin{defi}
The normalizer $N_C$ of  a subbundle $C$ of the omni-Lie algebroid $
\omni=\dev E\oplus \jet E$ is composed of all the sections of
$\omni$ that preserve $\Gamma(C)$ from the left side, that is,
\begin{equation}
N_C=\{X\in\Gamma(\omni)\mid \Dorfman{X,Y}\in\Gamma( C),\quad
\forall~Y\in\Gamma(C)\}.
\end{equation}
\end{defi}

It is easy to see that the  normalizer $N_C$ of $C$  is a Leibniz
subalgebra\footnote{Analogously, we may define  $N'_C$ to be the set
of sections of $\omni$ that preserve $C$ from the right side. But it
is not a Leibniz subalgebra.} of $\Gamma(\omni)$.

For any $X\in \Gamma(\omni)$,  we introduce the \emph{falling}
operator
$$
\falling{(\cdot)}:~~\ \Gamma(\omni)\longrightarrow \Gamma(\dev
\huaT),
$$
which is defined by
\begin{equation}
\falling{X}(t)\defbe \pomnib\Dorfman{X,Y},\quad \forall ~t\in
\Gamma(\mathcal{T}),
\end{equation}
where $Y\in \Gamma(\omni)$  satisfying $ \pomnib(Y)=t.$ By Lemma
\ref{lem: idea hom(TM,E)}, this    is well defined  and if $h\in
\Gamma(\Hom(\huaT,E))$,  $\falling{h}=h$.

In this section, we study the normalizer $N_{L}$ of a Dirac
structure $L$. Using the falling operator defined above, we
establish the relation between the normalizer $N_{L}$ of a reducible
Dirac structure $L$ and the derivation algebra $\Der(\pomnib (L))$
of the projective Lie algebroid $\pomnib (L)$.

\begin{pro}\label{pro:property of falling}
The falling operator $\falling{(\cdot)}$ is a morphism of Leibniz
algebras. Furthermore, $\forall$ $X\in\Gamma(\omni),
~t\in\Gamma(\huaT)$, we have
\begin{equation}\label{eqt:property falling}
pr_{TM}(\falling{X}(t))=[\jd\circ\rho
(X),pr_{TM}(t)]=[\jd(\falling{X}),pr_{TM}(t)].
\end{equation}
Conversely, given any $\delta\in \Gamma(\dev \huaT)$ satisfying
Equation \eqref{eqt:property falling}, there exists an
$X_\delta\in\Gamma(\omni)$ such that $\falling{X_\delta}=\delta$.
\end{pro}
\pf For all $X,Y\in \Gamma(\omni),~t\in\Gamma(\huaT)$ and
$Z\in\Gamma(\omni)$ satisfying $\pomnib(Z)=t$,  we have
\begin{eqnarray*}
\falling{\Dorfman{X,Y}}(t)&=&\pomnib\Dorfman{\Dorfman{X,Y},Z}=\pomnib(\Dorfman{X,\Dorfman{Y,Z}}-\Dorfman{Y,\Dorfman{X,Z}}\\
&=&\falling{X}\pomnib\Dorfman{Y,Z}-\falling{Y}\pomnib\Dorfman{X,Z}\\
&=&\falling{X}\circ\falling{Y}(t)-\falling{Y}\circ\falling{X}(t)\\
&=&[\falling{X},\falling{Y}]_{\dev }(t),
\end{eqnarray*}
which implies that the falling operator $\falling{(\cdot)}$ is a
morphism of Leibniz algebras.

Given  $X=\frkd+\mu$ and $Z=\frkr+\nu$ such that
$pr_{TM}(t)=\jd(\frkr)$, we have
\begin{eqnarray*}
pr_{TM}(\falling{X}(t))&=&pr_{TM}\pomnib\Dorfman{X,Z}\\
&=&[\jd(\frkd),\jd(\frkr)]=[\jd\circ\rho(X),pr_{TM}(t)]\\
&=&[\jd(\falling{X}),pr_{TM}(t)],
\end{eqnarray*}
which implies  Equation (\ref{eqt:property falling}).

Suppose that $\delta\in \Gamma(\dev \huaT)$ satisfies  Equation
(\ref{eqt:property falling}). Write $x=\jd(\delta)$   and define
$\chi=pr_E\circ\delta$. One has $$ \chi(ft)=x(f)
pr_E(t)+f\chi(t),\quad \forall ~f\in C^\infty(M).
$$
Therefore, $\chi|_{\frkX(M)}$ is $\CWM$-linear and there is an
associated $X_M\in \Gamma(\Hom(TM,E))$. Moreover, $\chi|_{\Gamma(
E)}$ is a derivation and there is an associated  $X_E\in\Gamma(\dev
E)$ such that $\jd(X_E)=x$. In turn, the operation of $\delta$ can
be expressed as
$$
\delta(y+v)=[x,y]+X_E(v)+X_M(y),\quad\forall~ y+v\in \Gamma(\huaT).
$$
Let $X_\delta= X_E+X_M \in \Gamma(\dev
E)\oplus\Gamma(\Hom(TM,E))\subset
 \Gamma(\omni)$. We claim that $\falling{X_\delta}=\delta$. In fact,
for any $ y+v\in \Gamma(\mathcal{T})=\Gamma(TM)\oplus\Gamma(E)$ and
$Y=\frkr+\nu\in \Gamma(\omni)$ satisfying $\jd(\frkr)=y$ and
$\p(\nu)=v $, we have
\begin{eqnarray*}
 \falling{X_\delta}(y+v)&=&\pomnib\Dorfman{X,Y}=\pomnib([X_E,\frkr]_\dev+\LieDerivation_{X_E}\nu
 -\LieDerivation_{\frkr}X_M + \jetd(X_M(y)))\\
 &=&[x,y]+X_E(v)+X_M(y) =\delta(y+v).
\end{eqnarray*}
\hfill$\square$

Let $A\subset\huaT$ be a  projective Lie algebroid and $\Inn (A)$
the set of inner derivations, which  consists of operators
$[a,\cdot]_A,$ where $a\in \Gamma(A)$. Denote the set of external
derivations by $\Ext(A)$, i.e.
\begin{equation}\label{eqt:out A}
\Ext(A)=\Der(A)/\Inn (A).
\end{equation}
By Theorem \ref{Thm:Main}, there is a unique lifted Dirac structure
$L^A$ such that $A$ is the quotient Lie algebroid of $L^A$.
Concerning the relation between the normalizer $N_{L^A}$ and  the
derivation algebra $\Der(A)$, we have

\begin{thm}\label{thm:Normal}
If $X\in N_{L^A}$, then $\falling{X}\mid_A\in \Der(A)$. Conversely,
for any $\delta\in\Der(A)$, there exists an $X_\delta\in N_{L^A}$,
such that $\falling{(X_\delta)}\mid_A=\delta$. Moreover, we have the
following commutative diagram where the two rows are exact
sequences:
$$
\xymatrix{
 0 \ar[r] & \Gamma(A^0)  \ar[d]_{i} \ar[rr]^{i} &&
                 \Gamma(L^A)  \ar[d]_{{{i}}}  \ar[rr]^{\falling{(\cdot)}\mid_A} && \Inn(A) \ar[d]_{{i}} \ar[r]  & 0\\
 0 \ar[r] & \Gamma(A^0)\oplus\Gamma(E)   \ar[rr]^{\kappa} &&
                 N_{L^A}  \ar[rr]^{\falling{(\cdot)}\mid_A } && \Der(A)
                  \ar[r]  & 0.
                }
$$
Here $i$ is the inclusion. The map $\kappa$ is defined by
$\kappa(\phi+u)=\phi+\jetd u$, $\forall$ $\phi\in
\Gamma(A^0),~u\in\Gamma(E)$. In particular, $\falling{X}\mid_A\in
\Inn A$ if and only if $X=l+\jetd u$, for some $l\in \Gamma(L^A)$,
$u\in \Gamma(E)$.
\end{thm}
\pf If $X\in N_{L^A}$,  then for any $a_1$, $a_2\in \Gamma(A)$, we
can find $l_1,l_2\in \Gamma(L^A)$ such that $\pomnib(l_i)=a_i$.
Hence
\begin{eqnarray*}
\falling{X}[a_1,a_2]_A&=& \pomnib\Dorfman{X,\Dorfman{l_1,l_2}}
=\pomnib\Dorfman{\Dorfman{X,l_1},l_2}
+\pomnib\Dorfman{l_1,\Dorfman{X,l_2}}\\
&=&[\pomnib\Dorfman{X,l_1},a_2]_A+[a_1,\pomnib\Dorfman{X,l_2} ]_A\\
&=&[\falling{X}a_1,a_2]_A+[a_1,\falling{X}a_2]_A\,,
\end{eqnarray*}
which implies that $\falling{X}\mid_A\in \Der(A)$.

Conversely, given any $\delta\in \Der(A)$, set $x=\jd(\delta)\in
\frkX(M)$ and find an extension  $\widetilde{\delta}\in\Der(\huaT)$
of $\delta$, that is,   $\alpha(\widetilde{\delta})=x$ and
$\widetilde{\delta} |_{\Gamma(A)}=\delta$. Since the elements of
$\Der(A)$  satisfy (\ref{eqt:property falling}) and by Proposition
\ref{pro:property of falling}, there is an
$X_{\widetilde{\delta}}=X_E+X_M$ such that
$\falling{X_{\widetilde{\delta}}}=\widetilde{\delta}$, i.e.
$\falling{X_{\widetilde{\delta}}}\mid_A=\delta$.

Next we prove $X_{\widetilde{\delta}}\in N_{L^A}$. For all $l\in
\Gamma(L^A)$, it is evident that
$\Dorfman{X_{\widetilde{\delta}},l}\in\Gamma(\pomnib\inverse(A))$.
Furthermore, $\forall$ $Y\in \Gamma(\pomnib\inverse(A))$, we have
\begin{eqnarray*}
\pomnib\Dorfman{\Dorfman{X_{\widetilde{\delta}},l},Y}&=&
\pomnib\Dorfman{X_{\widetilde{\delta}},\Dorfman{l,Y}}-\pomnib\Dorfman{l,\Dorfman{X_{\widetilde{\delta},Y}}}\\
&=&\falling{X_{\widetilde{\delta}}}[\pomnib l,\pomnib Y]_A-[\pomnib
l,
\falling{X_{\widetilde{\delta}}}(\pomnib Y)]_A\\
&=&[\falling{X_{\widetilde{\delta}}}(\pomnib l),\pomnib(Y)]_A
=[\pomnib\Dorfman{X_{\widetilde{\delta}},l},\pomnib Y]_A,
\end{eqnarray*}
which implies that $X_{\widetilde{\delta}}\in N_{L^A}$.

For an $X\in N_{L^A}$ satisfying $\falling{X}(\Gamma(A))=0$,  it is
easy to see that $\jd(\falling{X})=0$, i.e. $\jd\circ\rho(X)=0$. So
we are able to write
$$
X=\Phi+\frky+\jetd u,\quad\mbox{ where } \Phi\in \Gamma(\gl(E)),~
\frky\in \Gamma(\Hom(TM,E)),~ u\in \Gamma(E).
$$
Clearly, $\Dorfman{\jetd u,~\cdot~}=0$. By Lemma \ref{lem: idea
hom(TM,E)}, we have $ \Phi+\frky \in \Gamma(A^0)$, which implies
that $\ker(\falling{X}\mid_A)=\Gamma(A^0)\oplus\Gamma(E)$. The
remaining statements of the theorem are easy to be checked and we
omit the details. \hfill$\square$


\begin{ex}\rm For a reducible Dirac structure $L_{\pi}$ given in Example
\ref{Ex:Graphpi}, we consider its normalizer. For $ u\in\Gamma(E)$,
since we have $\{\jetd u,L_{\pi}\}=0$,  it suffices to consider
elements of the form $\frkd + \frky\in \Gamma(\omni)$, where
$\frky\in \Gamma(\Hom(TM,E))$. Rewrite $\frkd
+\frky=\frkd-\pi(\frky)+\pi(\frky)+\frky$ where $\pi(\frky)+\frky\in
L_{\pi}$. We have
$$
\frkd+\frky\in N_{L_{\pi}}\Longleftrightarrow\frkd-\pi(\frky)\in
N_{L_{\pi}}\Longleftrightarrow\LieDerivation_{\frkd-\pi(\frky)}\circ\pi=\pi\circ\LieDerivation_{\frkd-\pi(\frky)}
\Longleftrightarrow d(\frkd-\pi(\frky))=0.
$$

Here the coboundary operator $d$ is associated with cochain complex
$\Gamma(\Hom(\wedge^\bullet\jet E,E))$ and the representation
$\pi:~\jet E\longrightarrow\dev E$, which is known as the adjoint
representation \cite{CFSecondary,Shengdeformation}. From this we get
the following exact sequence:
$$
0 \rightarrow \Gamma(\Hom(TM,E))\oplus\Gamma(E)
\stackrel{\kappa}{\longrightarrow}
                  N_{L_{\pi}} \stackrel{p}{\longrightarrow} B(\jet E, E)\cap\dev
E \rightarrow
                  0,
$$
where $B(\jet E, E)$ is the set of $1$-cocycles and the maps
$\kappa,~p$ are given by
$$\kappa(\frky+u)=\pi(\frky)+\frky+\jetd u,\quad
p(\frkd+\frky+[u])=\frkd-\pi(\frky),$$ where
$\frky\in\Gamma(\Hom(TM,E)),~\frkd\in\Gamma(\dev E),~u\in\Gamma(E)$.
\end{ex}


\section{Cohomology of Dirac Structures}\label{sec:2cohomology}

By Proposition \ref{pro:representation}, for any Dirac structure
${L}\subset \omni$, there is a representation $\rho_{L}$ on $E$. Let
$d_{L}:\Gamma(\Hom(\wedge^\bullet {L},E))\lon
\Gamma(\Hom(\wedge^{\bullet+1} {L},E))$ be the associated coboundary
operator. In this section, we study  the cohomology group
$\mathrm{H}^\bullet(L,\rho_{L})$ and explore the relation between
$N_{L}$ and $\mathrm{H}^1(L,\rho_{L})$. We also study the
deformation  of a Dirac structure, which is related with
$\mathrm{H}^2(L,\rho_{L})$.

Let $L\subset\omni$ be a Dirac structure. For any
$X\in\Gamma(\omni)$, $
\omega_X=\ppairingE{X,~\cdot~}:~~{L}\longrightarrow E $ naturally
defines a $1$-cochain. We first prove the following fact.
\begin{pro}\label{pro:Normalizer}
    $X\in N_L\Longleftrightarrow
d_L\omega_X=0 $.
\end{pro}
\pf   For any $l_1, l_2\in \Gamma({L})$, we have
\begin{eqnarray*}
d_{L}\omega_X(l_1,l_2)&=&\rho_L(l_1)\omega_X(l_2)-\rho_L(l_2)\omega_X(l_1)-\omega_X(\Dorfman{l_1,l_2})\\
&=&\rho_L(l_1)\ppairingE{X,l_2}-\rho_L(l_2)\ppairingE{X,l_1}-\ppairingE{X,\Dorfman{l_1,l_2}} \\
&=&\ppairingE{\Dorfman{l_1,X},l_2}-\ppairingE{\Dorfman{l_2,X},l_1}-\ppairingE{X,\Dorfman{l_2,l_1}}\\
&=&-\ppairingE{\Dorfman{X,l_1},l_2}+\ppairingE{\Dorfman{X,l_2},l_1}-\ppairingE{X,\Dorfman{l_2,l_1}}\\
&&+\ppairingE{2\jetd \ppairingE{l_1,X},l_2}-\ppairingE{2\jetd
\ppairingE{l_2,X},l_1}
\\
&=&-2\ppairingE{\Dorfman{X,l_1},l_2}
-\ppairingE{X,\Dorfman{l_2,l_1}}
-\rho_L(l_1)\ppairingE{X,l_2}+\rho_L(l_2)\ppairingE{X,l_1}
\\
&=&-2\ppairingE{\Dorfman{X,l_1},l_2}-d_L\omega_X(l_1,l_2).\quad\quad
\end{eqnarray*}
Therefore,
\begin{equation}\label{Eqt:domegaX}
d_{L}\omega_X(l_1,l_2)=-\ppairingE{\Dorfman{X,l_1},l_2}.
\end{equation}
Since $L^\bot=L$, the above equality implies that $X\in
N_L\Longleftrightarrow d_L\omega_X=0.$ \hfill$\square$

\begin{pro}   Let $A\subset\huaT$ be a projective Lie algebroid and $L^A $  the lifted Dirac
structure, for any $X\in N_{L^A}$, $\omega_X$ is a coboundary if and
only if
 $
\falling{X}\in \Inn A$
\end{pro}
\pf By definition, $ \omega_X=d_{L^A} u$, for some $u\in\Gamma(E)$,
if and only if
$$
\ppairingE{X-2\jetd u,L^A}=0\quad\Longleftrightarrow X=2\jetd u +l,
\quad\mbox{ for some } l\in \Gamma(L^A).
$$
So the conclusion follows directly by Theorem  \ref{thm:Normal}.
\hfill$\square$

\begin{cor}
With the above notations, there is a natural inclusion
$$i:\Ext(A)\longrightarrow \mathrm{H}^1(L^A,\rho_{L^A}),
~~~~~~~~~~~\, \quad ~~~~~ i(\delta)=\omega_{X_\delta}, ~~~~~~~~
\forall~ \delta\in\Der(A).
$$
where $\Ext(A)$ is defined by \eqref{eqt:out A} and $X_\delta $ is
given in Theorem \ref{thm:Normal}.
\end{cor}
\pf By Theorem \ref{thm:Normal}, for any $\delta\in\Der(A)$, there
is an $X_\delta\in N_{L^A}$ such that
$\falling{X_\delta}|_A=\delta$. By Proposition \ref{pro:Normalizer},
  $ i(\delta)=\omega_{X_\delta} $ is closed.


To see that $i$ is well defined, we  note that $\omega_{X_\delta}$
does not depend on the choice of $X_{\delta}$ (by Theorem
\ref{thm:Normal}). And if $\delta\in \Inn
 (A)$, then by Theorem \ref{thm:Normal} again,  $X_\delta=l+2\jetd u$, where $l\in
 \Gamma(L^A)$ and $u\in\Gamma(E)$. Therefore, $\omega_{X_\delta}=d_{L^A}u$ is exact.

Finally, the previous proposition implies that $i$ is injective.
\hfill$\square$
 \vspace{3mm}

Suppose that $E$ and $\EStar$ are both Lie algebroids, respectively,
with anchors
 $\alpha$ and $\alpha^*$. Let $d_*:\Gamma(\wedge^\bullet E)\lon
\Gamma(\wedge^{\bullet+1} E)$ be the  Lie algebroid coboundary
operator associated with the Lie algebroid structure on $\EStar$. So
we have $d_*^2=0$. By definition, $(E,\EStar)$ is a Lie bialgebroid
if the following equality holds:
\begin{equation}\label{eqt:Lie bi}
d_*[u,v]=[d_*u,v]+[u,d_*v],\quad \forall~u,~v\in\Gamma(E).
\end{equation}
  (For more details about  Lie
bialgebroids, see \cite{Mkz:GTGA} and \cite{LWXmani}). The operator
$d_*:~\Gamma(E)\lon \Gamma(\wedge^{2} E)$ can be lifted to a bundle
map $\hat{d_*}:~\jet E\lon \wedge^2 E$, defined
 by
\begin{equation}\label{eqt:d hat}
\hat{d_*}(\jetd u)\defbe d_* u,\ \ \hat{d_*}(\dM f\otimes u)\defbe
{d_*}f\wedge u,\quad\forall ~u\in \Gamma(E),~ f\in \CWM.
\end{equation}
In \cite{CLomni}, we  proved that a Lie algebroid structure on $E$
can be lifted to a bundle map $\pisharp:~\jet E\lon \dev E$, which
is also a representation of the jet Lie algebroid $(\jet
E,\pibracket{\cdot,\cdot},\jd\circ\pisharp)$ on $E$, where $\pi$ is
given by $\pi(\jetd u)(v)=[u,v]$ (known as the adjoint
representation of a Lie algebroid) and the Lie bracket
$\pibracket{\cdot,\cdot}$ is given by (\ref{pibracket}). So we have
an induced tensor representation $\widetilde{\pi}$ of $\jet E$ on
$\wedge^2 E$ given by
$$
 \widetilde{\pi}(\jetd u)(\huaW) = [u,\huaW ],\quad
\widetilde{\pi}(\dM f\otimes u)(\huaW) = [\huaW,f]\wedge
u,\quad\forall ~\huaW\in\Gamma(\wedge^2 E).
$$

\begin{pro}~~
\begin{itemize}
\item[1)] The pair $(E,\EStar)$ is a Lie bialgebroid   if and only if
$\hat{d_*}$ is a 1-cocycle.
\item[2)] The pair $(E,\EStar)$ is a coboundary
Lie bialgebroid (i.e. $d_*=[\tau,~\cdot~]$, for some $\tau\in
\Gamma(\wedge^2 E)$) if and only if $\hat{d_*}$ is a coboundary.
\end{itemize}
\end{pro}
\pf For all $u,~v\in \Gamma(E)$ and $f,~g\in C^\infty(M)$, we have
the following three formulas which  are given in \cite{CLomni}:
\begin{eqnarray*}
\pibracket{\jetd u,\jetd v}&=&\jetd[u,v],\\
\pibracket{\jetd u,\dM f\otimes v}&=& \dM\rho(u)(f) \otimes v+ \dM
f\otimes
[u,v],\\
\pibracket{\dM f\otimes u, \dM g\otimes v}&=& \rho(u)(g) (\dM
f\otimes v) -\rho(v)(f) (\dM g\otimes u).
\end{eqnarray*}

Denote the coboundary operator associated with the representation
$\widetilde{\pi}$ by $\huaD$. We have
\begin{eqnarray*} \huaD(\hat{d_*})(\jetd u,\jetd
v)&=&\widetilde{\pi}(\jetd u)\hat{d_*}(\jetd
v)-\widetilde{\pi}(\jetd v)\hat{d_*}(\jetd u)-\hat{d_*}([\jetd
u,\jetd v]_\pi)\\
&=&[u,d_*v]-[v,d_*u]-d_*[u,v],\\
 \huaD(\hat{d_*})(\jetd u,\dM f\otimes v)
&=&(d_*[u,f]-[d_*u,f]-[u,d_*f])\wedge v,\\
 \huaD(\hat{d_*})(\dM f\otimes u,\dM g\otimes v)
&=&([d_*f,g]+[f,d_*g])\wedge u\wedge v.
\end{eqnarray*}
which implies that $(E,\EStar)$ is a Lie bialgebroid   if and only
if $\hat{d_*}$ is closed.

It is clear that $d_*=[\tau,~\cdot~]\Longleftrightarrow
\hat{d_*}=\huaD\tau$, which implies   (2). \hfill$\square$
\vspace{3mm}

Finally we consider the deformation of a projective Lie algebroid
$A$ and its lifted Dirac structure $L^A$. Let $\Omega:A\wedge
A\longrightarrow A\cap E$ be a bundle map. Consider an
$\varepsilon$-parameterized family of brackets
$$
[a,b]_A^\varepsilon=[a,b]_A+\varepsilon\Omega(a,b),\quad\forall~a,b\in
\Gamma(A).
$$
If every $\varepsilon$- bracket  endows $A$ a projective Lie
algebroid structure, we say that $\Omega$ generates a deformation of
the projective Lie algebroid $A$. Evidently, this requirement is
equivalent to the following compatibility conditions:
\begin{eqnarray}
\label{eqt:Omega
closed}\Omega([a,b]_A,c)+[\Omega(a,b),c]_A+~c.p.~&=&0,\\\label{eqt:Omega
fibrewise} \Omega(\Omega(a,b),c)+~c.p.~&=&0.
\end{eqnarray}
 Equation (\ref{eqt:Omega fibrewise}) means that $\Omega$ itself
defines a (fibrewise) Lie bracket. Furthermore, for all
$l_1,~l_2,~l_3\in \Gamma(L^A)$, we have
\begin{eqnarray*}
d_{L^A}\pomnib^*\Omega(l_1,l_2,l_3)&=&\rho_{L^A}(l_1)\pomnib^*\Omega(l_2,l_3)+c.p.+\pomnib^*\Omega(\Dorfman{l_1,l_2},l_3)+c.p.\\
&=&\rho_{L^A}(l_1)\Omega(\pomnib l_2,\pomnib
l_3)+c.p.+\Omega(\pomnib\Dorfman{l_1,l_2},\pomnib l_3)+c.p.\\
&=&\pomnib\Dorfman{l_1,\jetd \Omega(\pomnib l_2,\pomnib
l_3)}+c.p.+\Omega([\pomnib l_1,\pomnib l_2]_A,\pomnib l_3)+c.p.\\
&=&[\pomnib l_1,\Omega(\pomnib l_2,\pomnib
l_3)]_A+c.p.+\Omega([\pomnib l_1,\pomnib l_2]_A,\pomnib l_3)+c.p.,
\end{eqnarray*}
which implies that Equation (\ref{eqt:Omega closed}) is equivalent
to the requirement that  $\pomnib^*\Omega$ is closed.

Since there is a one-to-one correspondence between  reducible Dirac
structures
 and  projective Lie algebroids, we can associate a
deformation of the Dirac structure $L^A$ to the deformation of the
projective Lie algebroid $A$. Denote the deformed projective Lie
algebroid by $A_\varepsilon$, then the deformed Dirac structure
$L^{A_\varepsilon}$ is give by
$$
L^{A_\varepsilon}=\set{l+h~|~ l\in L_A, h \in \Hom(\huaT,E), \mbox{
s.t., } h(a)=\Omega(\pomnib(l),a),\quad\forall~ a\in A }.
$$

An interesting problem is to consider a deformation $\Omega$ which
is a coboundary:
\begin{equation}\label{Eqt:temp314}\pomnib^*\Omega=d_{L^A}
\omega_X,\quad\mbox{ for some }X\in \Gamma(\omni).
\end{equation}
\begin{pro} Let
$\Omega:A\wedge A\longrightarrow A\cap E$ be a bundle map. If
$\pomnib^*\Omega=d_{L^A} \omega_X$ for some $X\in \Gamma(\omni)$,
then $X\in N_{A^0}=N_{\pomnib\inverse(A)}$. Moreover, we have
\begin{equation}\label{Eqt:lambdadomegaX}
\Omega(a_1,a_2)=\half([\falling{X}a_1,a_2]_A+[a_1,\falling{X}a_2]_A-\falling{X}[a_1,a_2]_A)\,,\quad\forall~
a_1,a_2 \in \Gamma(A).
\end{equation}
Furthermore, $\Omega$ generates a deformation of the projective Lie
algebroid $A$ if and only if
\begin{equation}\label{Eqt:temp454}[T^X(a,b),c]_A+T^X([a,b]_A,c)+~c.p.~=0,\quad\forall ~a,b,c\in
\Gamma(A),\end{equation} where $T^X:~\Gamma(A)\wedge\Gamma(A)\lon
\Gamma(A)$ is defined by
$$
T^X(a,b)\defbe
\falling{X}([\falling{X}a,b]_A+[a,\falling{X}b]_A-\falling{X}[a,
b]_A)-[\falling{X}a,\falling{X}b]_A.
$$
Conversely, for any $X\in N_{\pomnib\inverse(A)}=N_{A^0}$ satisfying
\eqref{Eqt:temp454}, $\Omega$ defined by equation
\eqref{Eqt:lambdadomegaX} is a bundle map from $A\wedge A$ to $A\cap
E$ that generates a deformation of $A$ and relation
\eqref{Eqt:temp314} holds.
\end{pro}
\pf By Equations (\ref{Eqt:domegaX}) and (\ref{Eqt:temp314}), for
all $\theta\in \Gamma(A^0)$, $l\in \Gamma(L^A)$, we have
\begin{eqnarray*}
0=(\pomnib^*\Omega)(\theta,l)=d_{L^A}\omega_X(\theta,l)=
-\ppairingE{\Dorfman{X,\theta},l }.
\end{eqnarray*}
Thus, $\Dorfman{X,\theta}\in \Gamma(L^A \cap
\Hom(\huaT,E))=\Gamma(A^0)$, i.e. $X\in N_{A^0}$. For any $Y\in
\Gamma(\pomnib\inverse(A))$, we have
\begin{eqnarray*}
\Dorfman{X,\theta}(\pomnib Y)&=& 2\ppairingE{\Dorfman{X,\theta},Y}
= 2\rho(X)\ppairingE{\theta,Y}-2\ppairingE{\theta,\Dorfman{X,Y}}\\
&=&\rho(X)\theta(\pomnib Y)-\theta\circ\pomnib\Dorfman{X,Y},
\end{eqnarray*}
which implies that $\theta\circ\pomnib\Dorfman{X,Y}=0$, i.e. $X\in
N_{\pomnib\inverse(A)}$.

 Let $l_i\in\Gamma(L^A)$ and
$\pomnib(l_i)=a_i$. By some straightforward computation, we have
\begin{eqnarray*}
\Omega(a_1,a_2)&=&(\pomnib^*\Omega)(l_1,l_2)=d_{L^A}\omega_X(l_1,l_2)\\
&=&-\ppairingE{\Dorfman{X,l_1},l_2
}=-\half\pomnib(\Dorfman{\Dorfman{X,l_1},l_2}+
\Dorfman{l_2,\Dorfman{X,l_1}})\\
&=&
-\half\pomnib(\Dorfman{X,\Dorfman{l_1,l_2}}-\Dorfman{l_1,\Dorfman{X,l_2}}+
\Dorfman{l_2,\Dorfman{X,l_1}})\\
&=&\half([\falling{X}a_1,a_2]_A+[a_1,\falling{X}a_2]_A-\falling{X}[a_1,a_2]_A),
\end{eqnarray*}
which implies  Equation (\ref{Eqt:lambdadomegaX}).

If $\Omega $ generates a deformation of the projective Lie algebroid
$A$, $\Omega $ itself defines a fibrewise  Lie bracket. It is easy
to see that this is equivalent to (\ref{Eqt:temp454}). The other
conclusions can be easily checked. \hfill$\square$\vspace{3mm}

For a  Lie algebroid $(\huaA,[\cdot,\cdot],\alpha)$, a Nijenhuis
operator is a bundle map $N: ~\huaA\lon \huaA$ such that the
following equality holds
$$
T^N(a,b)\defbe N([Na,b]+[a,Nb]-N[a, b])-[Na,Nb]=0,\quad\forall ~a,b
\in \Gamma(\huaA).
$$
It induces a new Lie algebroid $(\huaA, [\cdot,\cdot]_N,\alpha_N)$,
where $\alpha_N=\alpha\circ N$ and
$$
[a, b]_N=[Na,b]+[a,Nb]-N[a, b].
$$
In fact, $T^N=0$ is only a sufficient condition for the bracket
operation $[\cdot,\cdot]_N$ being a Lie bracket. The necessary and
sufficient condition is
$$
[a,T^N(b,c)]+T^N(a,[b,c])+c.p.=0.
$$

 The role of the operator $\falling{X}:~ \Gamma(A)\lon \Gamma(A)$ is
just like that of  a Nijenhuis operator. In general, $\falling{X}$
is not a bundle map,  but it still induces a twist of the  Lie
algebroid. In fact, $\falling{X}$ is a bundle map if and only if
$X\in \Gamma(\Hom(\huaT,E))$ and in this case $\falling{X}=X|_A:~
A\lon A\cap E$, which is a Nijenhuis operator if $T^X$ vanishes.

\section*{Acknowledgments}
Z. Chen  would like to  thank   P. Xu and M. Stienon for
 the useful discussions and suggestions that helped him improving this work.
Y.-H. Sheng gives his warmest thanks to L. Hoevenaars, M. Crainic,
I. Moerdijk and C. Zhu for their help and useful comments during his
stay in Utrecht University and Courant Research Center,
G\"{o}ttingen, where a part of work was done. Research partially
supported by NSFC(10871007, 10911120391/A0109), Doctoral Fund. of
MEC (20090001110006) and CPSF(20060400017)(20090451267); the third
author financially supported by the governmental scholarship from
China Scholarship Council.

\end{document}